\documentclass[12pt]{article}

\usepackage{amsmath, amsfonts, amssymb}
\usepackage{latexsym}
\usepackage{graphicx}
\usepackage{color}
\usepackage{url, hyperref}
\usepackage{xcolor,soul,cancel}

\newtheorem{thm}{Theorem}
 \newtheorem{cor}[thm]{Corollary}
 \newtheorem{lem}[thm]{Lemma}
 \newtheorem{prop}[thm]{Proposition}
 \newtheorem{exa}[thm]{Example}

\newcommand{\al}{\alpha}
\newcommand{\be}{\beta}

\newcommand{\la}{\lambda}

\newcommand{\om}{\omega}

\newcommand{\si}{\sigma}
\newcommand{\Si}{\Sigma}

\newcommand{\n}{\|}
\newcommand{\vsk}{\vspace{1mm}}

\newcommand{\eq} [1] {\begin{equation}\label{#1}\quad}
\newcommand{\en} {\end{equation}}

\newcommand{\norm}[1]{\left\Vert#1\right\Vert}
\newcommand{\abs}[1]{\left\vert#1\right\vert}

\newcommand{\N}{\mathbb N}

\newcommand{\C}{\mathbb C}

\newcommand{\R}{\mathbb R}

\newcommand{\diag}{\operatorname{diag}}

\newcommand{\tr}{\operatorname{trace}}

\begin{document}

\begin{center}
{\LARGE \bf The norm attainment problem for functions of projections}

\vspace{5mm}
{\Large Albrecht B\"ottcher and Ilya M. Spitkovsky}

\vspace{8mm}
{\em \large For Bernd Silbermann on his 80th birthday}
\end{center}

\medskip
\begin{quote}
\footnotesize{The paper is concerned with the problem of identifying the norm attaining operators
in the von Neumann algebra generated by two orthogonal projections on a Hilbert space.
This algebra contains every skew projection on that Hilbert space and hence the results of the paper
also describe functions of skew projections and their adjoints that attain the norm.}


\end{quote}

This note is in the spirit of paper \cite{Sheck}.
The meta theorem of that paper is that the two projections
theorem of Halmos is something like Robert Sheckley's Answerer: no question
about the $W^*$- and $C^*$-algebras generated by two orthogonal projections
will go unanswered, provided the question is not foolish.
The norm attainment problem asks whether for a given bounded linear operator $A$ there is a unit vector $x$
such that $\n A x\n = \n A\n$. In this generality, a useful answer is not available.
Here we pose the question for the case where $A$ is a  function of two orthogonal projections or a function
of one skew projection and its adjoint.

\vsk
Let $P$ and $Q$ be orthogonal projections acting on a real or complex Hilbert space $\mathcal H$.
According to Halmos' ``Two projections theorem'' (see \cite{Halmos} and
consult \cite{BSpit10, Silb} for the history and  more on the subject),
there is a representation of $\mathcal H$ as an orthogonal sum
\eq{M} \mathcal H = M_{00}\oplus  M_{01}\oplus  M_{10}\oplus  M_{11}\oplus M\oplus  M \en
(the last two summands have the same dimension and are thus identified
via an appropriate unitary similarity) with respect to which
\eq{P} P = (1,1,0,0)\oplus \begin{bmatrix}I & 0 \\ 0 & 0\end{bmatrix} \en
and
\eq{Q} Q = (1,0,1,0)\oplus \begin{bmatrix}I-H & \sqrt{H(I-H)}\, \\  \sqrt{H(I-H)} & H\,\end{bmatrix}. \en
Here and below we use the string $(a_{00},a_{01},a_{00},a_{11})$ as an abbreviation for
$a_{00}I_{M_{00}}\oplus a_{01}I_{M_{01}}\oplus a_{10}I_{M_{10}}\oplus a_{11}I_{M_{11}}$,
while the blocks of the matrix component in \eqref{P},\eqref{Q} are operators on $M$.
The selfadjoint operator $H$ has the spectrum $\sigma(H)\subseteq [0,1]$, with $0,1$
not being its eigenvalues. In particular, $\min\si(H)<1$. Note that $M=\{0\}$ if and only if $P$ and $Q$ commute.
The other four subspaces $M_{ij}$ also may or may not be actually present in \eqref{M};
we will let $\Lambda=\{ (i,j)\colon \dim M_{ij}\neq 0\}$.

\vsk
According to the Giles-Kummer theorem (see \cite{GiKu} or \cite[Theorem 7.1]{BSpit10}),
the von Neumann algebra $W^*(P,Q)$
generated by $P$ and $Q$ consists of the operators $A$ admitting the representation
\eq{A} A = (a_{00},a_{01},a_{00},a_{11})\oplus
\begin{bmatrix}\phi_{00}(H) & \phi_{01}(H) \\ \phi_{10}(H) & \phi_{11}(H)\end{bmatrix} \en
with respect to the decomposition \eqref{M} of $\mathcal H$. Here $a_{ij}$ are
arbitrary complex numbers and $\phi_{ij}$ are (also arbitrary) functions in $L^\infty(\sigma(H))$ with
respect to the spectral measure of $H$. We will sometimes write the rightmost summand in \eqref{A} as $\Phi_A(H)$.

\vsk
The norms of operators from $W^*(P,Q)$ were computed in \cite[Theorem 10]{Spit94},
see also \cite[Theorem 7.9]{BSpit10}.
Namely, for $A$ as in \eqref{A},
\eq{norm} \norm{A}=\max\left\{\max_{(j,k)\in\Lambda}\abs{a_{jk}},\max_{x\in\sigma(H)}
\sqrt{\frac{\phi(x)+\sqrt{\phi(x)^2-4\abs{\omega(x)}^2}}{2}} \right\}. \en
Here $\phi=\sum_{i,j=0,1}\abs{\phi_{ij}}^2$ is the square of the Frobenius norm
of $\Phi_A$ and $\omega=\det\Phi_A=\phi_{00}\phi_{11}-\phi_{01}\phi_{10}$.
Note that always $\phi(x)^2 -4|\om(x)|^2 \ge 0$.

\vsk
The question we are addressing here is: when is $\norm{A}$ {\em attained}, i.e., when does
there exist a unit vector $x\in\mathcal H$ such that $\norm{Ax}=\norm{A}$? We will call $A$ a
{\em norm attaining} operator and write $A\in\mathcal N$ if this happens to be the case.

\vsk
For our purposes it is useful to recall that formula \eqref{norm} was derived in \cite{Spit94} from the fact that
\[ \lambda_{\max}:=\max_{x\in\sigma(H)}\frac{\phi(x)+\sqrt{\phi(x)^2-4\abs{\omega(x)}^2}}{2}\]
is the right endpoint of the spectrum $\sigma\left(\Phi_{A^*A}(H)\right)$.
Since for every operator $X$ acting on $\mathcal H$ we have $X\in\mathcal N$
if and only if $X^*X\in\mathcal N$ if and only if the right endpoint of $\sigma(X^*X)$ is its
eigenvalue (see \cite{CaNe} and \cite{PaPa}),
we just need to figure out when $\lambda_{\max}$ is (or is not) an eigenvalue of $A^*A$.

\vsk
To this end, recall that for operators in $W^*(P,Q)$ the description of their kernels is also
known (\cite[Theorem 7.5]{BSpit10} or \cite[Theorem 1]{Spit94}).
There is no need to include its exact form here, but an important for us consequence of it is as follows.

\begin{lem}\label{th:ker}The kernel of $\Phi_A(H)$ is non-trivial
if and only if the spectral measure of $\{x \in \si(H): \om(x)=0\}$ is non-zero.  \end{lem}

When applied to $A-\lambda I$ in place of $A$, Lemma~\ref{th:ker} immediately yields the following.

\begin{prop}\label{th:eig} Let $A$ be given by \eqref{A}. Then $\lambda\in\C$
is an eigenvalue of $A$ if and only if the spectral measure of the set
$\{x\in\sigma(H)\colon \lambda^2-\tr\Phi_A(x)\lambda+\det\Phi_A(x)=0\}$ is non-zero. \end{prop}

Since $\tr\Phi_{A^*A}=\phi$ and $\det\Phi_{A^*A}=\abs{\omega}^2$,
the eigenvalues $\lambda$ of $A^*A$ are characterized by the property that
the spectral measure of the set
$\{x\in\sigma(H)\colon \lambda^2-\phi(x)\lambda+\abs{\omega(x)}^2=0\}$ is non-zero.
In particular, $\lambda_{\max}$ is an eigenvalue of $A^*A$ (and not just a point of
its spectrum) if and only if the spectral measure of the set $x\in\sigma(H)$ on
 which the function \eq{psi} \psi(x):=\phi(x)+\sqrt{\phi(x)^2-4\abs{\omega(x)}^2} \en
 attains its maximum value is non-zero. Denoting this set by $\Sigma(A)$, we arrive at the following conclusion.

\begin{thm}\label{th:na} Let $A$ be the operator given by \eqref{A}. Then $A\in\mathcal{N}$ if and only if
either
{\em (i)} $\max_{(j,k)\in\Lambda}\abs{a_{jk}}\geq\sqrt{\lambda_{\max}}$
or
{\em (ii)} $\max_{(j,k)\in\Lambda}\abs{a_{jk}}<\sqrt{\lambda_{\max}}$ and the
spectral measure of $\Sigma(A)$ is non-zero. \end{thm}

Let now $T$ be a skew projection on $\mathcal H$. We assume that $T$ is genuinely skew, which
is equivalent to the requirement that $\n T\n >1$. Denote by $P=P_{{\rm Ran}\,T}$ the orthogonal projection onto
the range of $T$ and by $Q=P_{{\rm Ker}\,T}$ the orthogonal projection onto the kernel of $T$. Afriat \cite{Afriat}
(see also \cite[Proposition 1.6]{BSpit10}) showed that then $\n PQ\n < 1$ and
\begin{equation}T =(I-PQ)^{-1}P(I-PQ). \label{Af} \end{equation}
Moreover, in \eqref{M} then $M_{00}=M_{11}=0$, while $M_{01}$ and $M_{10}$
may or may not be present. From \eqref{P} and \eqref{Q} we obtain
\[ I-PQ= (1,1)\oplus\begin{bmatrix}H & -\sqrt{H(I-H)}\, \\ 0 & I \,\end{bmatrix},\]
and since $I-PQ$ is invertible, so also must be $H$. Formula \eqref{Af} then gives
\eq{T} T= (1,0)\oplus\begin{bmatrix}I & -\sqrt{H^{-1}(I-H)}\, \\ 0 & 0 \,\end{bmatrix}
=(1,0)\oplus\begin{bmatrix}I & -\sqrt{H^{-1}-I}\, \\ 0 & 0 \,\end{bmatrix}. \en

\begin{cor} \label{th:tna} A skew projection attains its norm if and
only if $\min\sigma(H)$ is an eigenvalue of the respective operator $H$.
\end{cor}

{\em Proof.} Indeed, for $A=T$ we get from \eqref{T} that $\phi(x)=x^{-1}$, $\omega(x)=0$, and hence $\psi(x)=2x^{-1}$,
which is a monotonically decreasing function. It follows that
$\Sigma(T)$ is the singleton $\{\min\sigma(H)\}$.
Since $\norm{T}>1$, part (i) of Theorem~\ref{th:na} is irrelevant
and the assertion follows from part (ii) of that theorem. $\;\:\square$

\vsk
From \eqref{P} and \eqref{Q} and the equalities $M_{00}=M_{11}=0$ we infer that
\[PQP=(0,0) \oplus \begin{bmatrix}I-H & 0\, \\ 0 & 0 \,\end{bmatrix}.\]
Thus, if $A=T$, we see that the respective operator $H$ in Corollary \ref{th:tna} is
\[H=(I-PQP)|{\rm Ran}\,P = I|{\rm Ran}\, T-P_{{\rm Ran}\,T}P_{{\rm Ker}\,T}|{\rm Ran}\,T.\]

We conclude with some examples.
The authors of \cite{Bala} recently proved that a skew projection $T$ is in $\mathcal{N}$ if and only
if the selfadjoint operator $T+T^*-I$ belongs to $\mathcal{N}$. The operator $T+T^*-I$ appeared in \cite{Buck}
and is therefore called the Buckholtz operator in \cite{Bala}. The following is an extension of this result.

\begin{exa}\label{Ex1} Let $T$ be a skew projection. Then the following are equivalent:

\smallskip
{\rm (i)} $T \in \mathcal{N}$,

\smallskip
{\rm (ii)} $T +\al T^*  +\be I \in \mathcal{N}$ for some $\al,\be \in \R$,

\smallskip
{\rm (iii)} $T +\al T^*  +\be I \in \mathcal{N}$ for all $\al,\be \in \R$.
\end{exa}

{\em Proof.} Fix $\al,\be \in \R$. We have to show that $T \in \mathcal{N} \Longleftrightarrow T+\al T^*+\be I \in \mathcal{N}$.
According to \eqref{T},
\[ A:=T+\al T^*+\be I= (1+\al+\be,\be)\oplus\begin{bmatrix}(1+\al+\be)I & -\sqrt{H^{-1}-I} \\-\al \sqrt{H^{-1}-I}  & \be I \end{bmatrix}. \]
Abbreviating $1+\al+\be$ to $s$ and $x^{-1}-1$ to $f(x)$ we obtain
\begin{eqnarray*}
& & \phi(x)=(1+\al+\be)^2+\be^2+(1+\al^2)(x^{-1}-1)=s^2+\be^2+(1+\al^2)f(x),\\
& & \om(x)=(1+\al+\be)\be -\al(x^{-1}-1)=s\be -\al f(x),\\
& & \psi(x)=s^2+\be^2+(1+\al^2)f(x)+\sqrt{(s^2+\be^2+(1+\al^2)f(x))^2-4(s\be-\al f(x))^2}.
\end{eqnarray*}
The term under the square root equals
\[
(s^2-\be^2)^2+2[(s^2+\be^2)(1+\al^2)+4s\be\al]f(x)+(1-\al^2)^2 f(x)^2.
\]
The function $f(x)$ is monotonically decreasing and nonnegative.
Since $2|s\be| \le s^2+\be^2$ and $2|\al|\le 1+\al^2$, we have $(s^2+\be^2)(1+\al^2)+4s\be\al \ge 0$.
Consequently, the term under the square root and therefore also $\psi(x)$ are monotonically decreasing,
which implies that $\Sigma(A)=\{\min\sigma(H)\}$.
Using that $\min \si(H) <1$, we get
\begin{eqnarray*}
\la_{\max} & \ge & \frac{s^2+\be^2+(1+\al^2)f(\min \si(H))+|s^2-\be^2|}{2}\\
& > &\frac{s^2+\be^2+|s^2-\be^2|}{2}=\max(s^2,\be^2),
\end{eqnarray*}
and Theorem  \ref{th:na}(ii) implies that $A \in \mathcal{N}$ if and only if the
spectral measure of the singleton
$\{\min \si(H)\}$ is positive, that is, if and only if $\min \si(H)$ is an eigenvalue of $H$.  This together with Corollary~\ref{th:tna}
shows that $A \in \mathcal{N} \Longleftrightarrow T \in \mathcal{N}$, as desired.
$\;\:\square$

\begin{exa} \label{Ex2}
Let $T$ be a skew projection.
Put $T^{(2)}=TT^*$, $T^{(3)}=TT^*T$, and more generally, $T^{(m)}:=\underbrace{TT^*TT^*T\cdots}_{m}$.
Then, for each $m$,  $T^{(m)} \in \mathcal{N} \Longleftrightarrow T \in \mathcal{N}$
\end{exa}

{\em Proof.} From \eqref{T} we infer that if $m=2k$ is even, then
\[T^{(m)}=(TT^*)^k= (1,0)\oplus\begin{bmatrix}H^{-k} & 0 \\ 0 & 0 \,\end{bmatrix}.\]
This implies that $\psi(x)=2x^{-2k}$ and thus $\Sigma(T^{(m)})=\{\min \si(H)\}$. Taking
into account that $\min \si(H) <1$, we obtain as above from Theorem \ref{th:na}(ii) and
Corollary~\ref{th:tna} that $T^{(m)} \in \mathcal{N} \Longleftrightarrow T\in \mathcal{N}$.
Finally, since $A\in \mathcal{N} \Longleftrightarrow AA^*\in \mathcal{N}$, we conclude that
\[T^{(2k+1)} \in \mathcal{N} \Longleftrightarrow T^{(4k+2)} \in \mathcal{N} \Longleftrightarrow T \in \mathcal{N},\]
which gives the assertion in the case of odd $m$. $\;\:\square$

\begin{exa} \label{Ex3} Let $\{\om_n\}_{n=1}^\infty$ be a sequence of positive real numbers that go monotonically to zero
and let $T$ be the skew projection on $\ell^2(\N)$ defined by the infinite matrix
\[T=\diag \left\{\begin{bmatrix}1 & -\om_n \\ 0  & 0 \end{bmatrix}_{n=1}^\infty\right\}.\]
Then $T \in \mathcal{N}$. Put $A=TT^*+T^*T-T-T^*-I$.  If $\om_n=1/n$, then $A \notin \mathcal{N}$,
but if $\om_n=2/n$, then $A \in \mathcal{N}$.
\end{exa}

{\em Proof.} It is clear that $T \in \mathcal{N}$: the norm is attained at the vector
\[(1/\sqrt{1+\om_1^2}, -\om_1/\sqrt{1+\om_1^2}, 0,0, \ldots)^\top \in \ell^2(\N).\]
Similarly one can treat the operator $A$ by sole calculations with $2\times 2$ matrices.
Here is how Theorem \ref{th:na} works. We may write
\[T=\begin{bmatrix}I & -\sqrt{H^{-1}-I}\, \\ 0 & 0 \,\end{bmatrix} \quad\mbox{with}\quad
H=\diag(x_n)_{n=1}^\infty, \]
where $\sqrt{x_n^{-1}-1}=\om_n$, that is, $x_n=1/(1+\om_n^2)$. Straightforward computation gives
\[A=\begin{bmatrix}H^{-1}-2I & 0\, \\ 0 & H^{-1}-2I \,\end{bmatrix}.\]
Thus, $\psi(x)=2(2-x^{-1})^2$. If $\om_n=1/n$, then
\[\si(H)=\left\{\frac{n^2}{n^2+1}: n \in \N\right\} \cup \{1\} =\left\{\frac{1}{2}, \frac{4}{5},\frac{9}{10}, \ldots\right\} \cup \{1\}.\]
The function $\psi(x)$ takes its maximum on $\si(H)$ at $x=1$, with $\la_{\max}=\psi(1)/2=1$.
Hence $\Si(A)=\{1\}$, and as $1$ is not an eigenvalue
of $H$, Theorem \ref{th:na}(ii) implies that $A \notin \mathcal{N}$. If $\om_n=2/n$, we have
\[\si(H)=\left\{\frac{n^2}{n^2+4}: n \in \N\right\} \cup \{1\} =\left\{\frac{1}{5}, \frac{1}{2},\frac{9}{13}, \ldots\right\} \cup \{1\}.\]
This time $\psi(x)$ assumes its maximum at $x=1/5$, the value of maximum being $\la_{\max}=\psi(1/5)/2=3.24$.  It follows that $\Si(A)=\{1/5\}$,
and since $1/5$ is an eigenvalue of $H$, we deduce from Theorem \ref{th:na}(ii) that $A \in \mathcal{N}$.
$\;\:\square$

\vsk
The last example can be elaborated to great extent. However, we leave it with Israel M. Gelfand:
``Explain this to me on a simple example; the difficult example I will be able to do on my own.''
({\url{http://www.israelmgelfand.com/edu_work.html}})

\bigskip
A. B\"ottcher,
Fakult\"at f\"ur Mathematik,
TU Chemnitz,
09107 Chemnitz,
Germany

\smallskip
{\tt aboettch@mathematik.tu-chemnitz.de}

\bigskip
I. M. Spitkovsky,
Division of Science and Mathematics,
New York  University Abu Dhabi (NYUAD),
Saadiyat Island, P.O. Box 129188,
Abu Dhabi,
United Arab Emirates

\smallskip
{\tt ims2@nyu.edu, ilya@math.wm.edu, imspitkovsky@gmail.com}

\end{document}